\documentstyle[12pt]{article}

\input amssym.tex
\input amssym.def

\begin{document}

\def\bcw{\mathbin{\bigcirc\mkern-15mu\wedge}}


\def\r#1{{\mathop{#1}\limits^\circ}}
\def\ring{\r}


\def\nint{\mathbin{\int\mkern-18mu\diagup \;} }

\newcommand{\namelistlabel}[1] {\mbox{#1}\hfil}
\newenvironment{namelist}[1]{%
\begin{list}{}
{\let\makelabel\namelistlabel
\settowidth{\labelwidth}{#1}
\setlength{\leftmargin}{1.1\labelwidth}}
}{%
\end{list}}

\newtheorem{lem}{Lemma}
\newtheorem{thm}{Theorem}

\newtheorem{Lem}{Lemma}[section]
\newtheorem{Thm}{Theorem}[section]


\newtheorem{cor}{Corollary}[section]

\newtheorem{Cor}[Lem]{Corollary}


\newtheorem{Prop}[Lem]{Proposition}


\newtheorem{prop}{Proposition}[section]


\newtheorem{Pro}[Thm]{Proposition}



\newtheorem{co}[Lem]{Corollary}

\newtheorem{corr}[Thm]{Corollary}


\newtheorem{C}[Lem]{Claim}

\newtheorem{Def}{Definition}[section]
\newtheorem{remark}{Remark}[section]

\pagestyle{myheadings}
\markboth{\underline{\sc On the Collapse of Tubes Carried by 3D Incompressible Flows }}
{\underline{\sc On the Collapse of Tubes Carried by 3D Incompressible Flows}}

\parskip=10pt

\baselineskip=18pt
\thispagestyle{empty}
\begin{center}
{\Large \bf \sc  On the Collapse of Tubes Carried \\
by 3D Incompressible Flows} \\
{\sc by}\\
{\sc Diego Cordoba\footnote{This work was supported
initially by the American Institute of Mathematics}
\&
Charles Fefferman\footnote{Partially supported by NSF Grant DMS 0070692}
\addtocounter{footnote}{-2}\footnotemark }
\end{center}
\medskip
\setcounter{section}{-1}
\section{\sc Introduction}
The 3-dimensional incompressible Euler equation (``3D Euler'') is as follows:
\[
\begin{array}{lcl}
\left( \frac{\partial}{\partial t} + u \cdot \nabla_x \right)
\, u \, = \, -
\nabla_x p && \left( x \in {\Bbb{R}}^3 , \, t \geq 0 \right) \\
\\
\nabla_x \cdot u \, = \, 0 && ( x \in {\Bbb{R}}^3 , \, t \geq 0 ) \\
\\
u ( x , 0 ) \, = \, u^0 ( x ) && ( x \in {\Bbb{R}}^3 )

\end{array}
\]
\noindent
with $u^0$ a given, smooth, divergence-free, rapidly decreasing vector field
on ${\Bbb{R}}^3$.  Here, $u ( x , t )$ and $p ( x , t )$ are the
unknown velocity and pressure for an ideal,
incompressible fluid flow at zero viscosity.
An outstanding open problem is to determine whether a 3D Euler solution can develop a singularity at a finite time $T$.  A classic result of
Beale-Kato-Majda [1] asserts that,
if a singularity forms at time $T$, then the vorticity
$\omega ( x , t ) = \nabla_x \times u ( x , t )$
grows so rapidly that
\[
\displaystyle{\int\limits^{T}_{0}} \, \sup_{x} \, | \,
\omega ( x , t ) \, | \, dt \, = \, \infty.
\]
\noindent
In [2], Constantin-Fefferman-Majda showed that, if the velocity
remain bounded up to the time $T$ of singularity formation, then the
vorticity direction $\omega ( x , t ) \slash \, | \omega ( x , t ) |$
cannot remain uniformly Lipschitz continuous up to time $T$.

\medskip
One scenario for possible formation of a singularity in a
3D Euler solution is a
constricting vortex tube.  Recall that a {\it vortex line} in a fluid is an arc on an integral curve of the vorticity $\omega ( x , t )$ for fixed $t$, and a
{\it vortex tube} is a tubular neighborhood in ${\Bbb{R}}^3$ arising
as a union of vortex lines.  In numerical simulations of 3D Euler solutions, one routinely sees that
vortex tubes grow longer and thinner, while bending and twisting.  If
the thickness of a piece of a vortex tube becomes zero in finite time, then
one has a singular solution of 3D Euler.
It is not known whether this can happen.

\medskip
Our purpose here is to adapt our work [3], [4] on two-dimensional flows
to three dimensions, for application to 3D Euler.
We introduce below the notion of a ``regular tube''.
Under the mild assumption that
\[
\displaystyle{\int\limits^{T}_{0}} \sup_x \, | u ( x , t ) \, | d t \, < \, \infty ,\]
\noindent
we show that a regular tube cannot reach zero thickness at time $T$.
In particular, for 3D Euler solutions, a vortex tube
cannot reach zero thickness in finite time, unless it
bends and twists so violently that no part of it forms a regular tube.
This significantly sharpens the conclusion of [2] for possible singularities
of 3D Euler solutions arising from vortex tubes.  On the other hand,
[2] applies to arbitrary singularities of 3D Euler solutions,
while our results apply to ``regular tubes''.

\medskip
Although we are mainly interested in 3D Euler solutions, our result is
stated for arbitrary incompressible flows in 3 dimensions.  The
proof is simple and elementary.  The main novelty for readers familiar with
[3], [4] is that
we can adapt the ideas of [3] to three dimensions, even though
there is no scalar that plays the r\^{o}le of the stream function
on ${\Bbb{R}}^2$.
\section{ \sc Regular Tubes}
Let $Q = I_1 \times I_2 \times I_3 \subset {\Bbb{R}}^3$ be a closed
rectangular box (with $I_j$ a bounded interval), and let $T > 0$ be given.

\medskip
A {\it regular tube} is an open set $\Omega_t \subset Q$ parametrized
by time $t \in [ 0 , T )$, having the form
\begin{equation}
\Omega_t \, = \,
\left\{
\left(
x_1 , x_2 , x_3 \right)
\in Q \, : \,
\theta \left(
x_1 , x_2 , x_3 , t \right) \, < \, 0 \right\} \,
\end{equation}
\noindent
with
\begin{equation}
\theta \in C^1 ( Q \times [ 0 , T ) ),
\end{equation}
\noindent
and satisfying the following properties:
\begin{equation}
| \, \nabla_{x_1 , x_2} \, \theta \, | \ne 0 \ {\rm for} \
(x_1 , x_2 , x_3 , t ) \in Q \times [ 0 , T ) , \,
\theta ( x_1 , x_2 , x_3 , t ) \, = \, 0 \, ;
\end{equation}

\begin{equation}
\Omega_t ( x_3 ) : = \, \left\{
( x_1 , x_2 ) \in I_1 \times I_2 \, : \,
(x_1 , x_2 , x_3 ) \in \Omega_t \right\} \ {\rm is \ non-empty} ,
\end{equation}

\noindent
for all $x_3 \in I_3 , \, t \in [ 0 , T )$;
\begin{equation}
{\rm closure} \ \left( \Omega_t ( x_3 ) \right) \subset \ {\rm interior} \
(I_1 \times I_2 )
\end{equation}
\noindent
for all $x_3 \in I_3$, $t \in [ 0 , T)$.

\medskip
For example, a thin tubular neighborhood of a curve $\Gamma$ forms a
regular tube, provided the tangent vector $\Gamma^\prime$ stays
transverse to the $( x_1 , x_2 )$ plane.

\medskip
Let $u ( x , t ) =
( u_k ( x , t ))_{1 \leq k \leq 3}$ be a $C^1$ velocity field
defined on $Q \times [ 0 , T )$.  We say that the regular tube
$\Omega_t$ {\it moves with the velocity field $u$},
if we have
\begin{equation}
\left(
\frac{\partial}{\partial t} \, + \,
u \cdot \nabla_x \right) \theta = 0 \ {\rm whenever} \
( x , t ) \in Q \times [ 0 , T ) ,
\theta ( x , t ) = 0 \, .
\end{equation}

\noindent
It is well-known that a vortex tube arising from a 3D Euler solution
moves with the fluid velocity.
\section{\sc Statement of the Main Result}

\medskip
\noindent{\bf Theorem:} \hspace{.5em}
{\it Let $\Omega_t \subset Q ( t \in [ 0 , T ))$ be a regular tube
that moves with a $C^1$, divergence free velocity field $u ( x , t )$}.

\medskip
\noindent
If
\begin{equation}
\displaystyle{\int\limits^{T}_{0}} \sup_{x \in Q} \, | u ( x , t ) \, | d t < \infty
\end{equation}

\noindent
then
\begin{equation}
\liminf_{t \to T-} \ {\rm Vol} (\Omega_t) \, > \, 0.
\end{equation}
\section{\sc Calculus Formulas for Regular Tubes}
Let $\Omega_t$ be a regular tube, as in (1)$\cdots$(5).  Recall that
\begin{equation}
\Omega_t ( x_3 ) = \left\{
(x_1 , x_2 ) \in I_1 \times I_2 : \, \theta ( x_1 , x_2 , x_3 , t ) \, < \, 0 \right\}.
\end{equation}

\noindent
Define also
\begin{equation}
S_t ( x_3) = \{ (x_1 , x_2 ) \in \ {\rm interior} \ (I_1 \times I_2):
\, \theta (x_1 , x_2 , x_3 , t ) = 0 \} \
{\rm for} \
x_3 \in I_3 , \, t \in [ 0 , T ) .
\end{equation}

\noindent
Also, for intervals $I \subset I_3$, and for $t \in [ 0 , T )$, define
\begin{equation}
\Omega_t ( I ) = \{
( x_1 , x_2 , x_3 ) \in Q : x_3 \in I \ {\rm and } \
\theta ( x_1 , x_2 , x_3 , t ) \, < \, 0 \} , \ {\rm and}
\end{equation}
\begin{equation}
S_t ( I ) = \{
(x_1 , x_2 , x_3 ) \in Q : \, x_3 \in I \ {\rm and} \
(x_1 , x_2 ) \in S_t (x_3) \} .
\end{equation}

\noindent
Let $\nu$ denote the outward-pointing unit normal to $S_t ( I_3)$,
and let $\tilde{\nu} = ( \tilde{\nu}_1 , \tilde{\nu}_2 , 0 )$,
where $(\tilde{\nu}_1 , \tilde{\nu}_2 )$ is the outward-pointing
unit normal to $S_t ( x_3)$.  Thus, $\nu$ and
$\tilde{\nu}$ are continuous vector-valued functions, defined
on ${\cal{S}} = \{
( x_1 , x_2 , x_3 , t ) \in Q \times [ 0 , T ) : ( x_1 , x_2 , x_3 ) \in S_t ( I_3 ) \}$.
Define also scalar-valued functions $\sigma , \tilde{\sigma}$ on
${\cal{S}}$ by requiring that
\begin{equation}
\left(
\frac{\partial}{\partial t} \, + \,
\sigma \nu \cdot \nabla_x \right) \, \theta =
\left(
\frac{\partial}{\partial t} \, + \, \tilde{\sigma} \tilde{\nu} \cdot \nabla_x \right) \, \theta \, = \, 0 \ {\rm for } \ x \in S_t ( I_3) .
\end{equation}

\noindent
Again, $\sigma$ and $\tilde{\sigma}$ are well-defined and continuous
on ${\cal{S}}$, thanks to (3).  Let $F$ be any continuous function on $Q$.  We will establish the following elementary formulas.
\begin{equation}
\frac{d}{dt} \left[ \:
\displaystyle{\int\limits_{\Omega_t (x_3)}}
 \,
F ( x_1 , x_2 , x_3 ) \, d \,  {(\rm Area)}
\right] \, = \,
\displaystyle{\int\limits_{(x_1 , x_2 ) \in S_t (x_3)}} \;
F(x_1 , x_2 , x_3 ) \,
\tilde{\sigma} ( x_1 , x_2 , x_3 , t) \, d \,
{\rm ( length)} \\
\end{equation}

\begin{equation}
\displaystyle{\int\limits_{S_t (I)}} \,
F d {\rm (Area)} \, = \,
\displaystyle{\int\limits_{x_3 \in I}} \,
\left\{ \:
\displaystyle{\int\limits_{S_t ( x_3)}} \:
\frac{F}{\nu \cdot \tilde{\nu}} \, d \,
{\rm (length)} \right\} \, d x_3 \, .
\end{equation}

\noindent
To check (14) and (15), we may use a partition of unity to reduce
to the case in which $F$ is supported in a small neighborhood $U$.  Also,
we may restrict attention to a small time interval $J$.  In a small enough
$U \times J$, we may assume that $S_t (I_3)$ is given by the graph of a
$C^1$ function $\psi$, thanks to (3).  Thus, without loss of generality,
we may suppose that, in $U \times J$, we have

\begin{equation}
(x_1 , x_2 , x_3 ) \in \Omega_t (I_3) \
{\rm  if \  and \  only \  if } \ x_1 < \psi ( x_2 , x_3 , t ) , \
{\rm and}
\end{equation}
\begin{equation}
(x_1 , x_2 , x_3 ) \in S_t (I_3) \ {\rm if \ and \ only \ if} \
x_1 \, = \, \psi ( x_2 , x_3, t ) .
\end{equation}

\noindent
From (16), (17) we have also
\begin{equation}
(x_1 , x_2 ) \in \Omega_t (x_3) \ {\rm if \ and \ only \ if} \
x_1 < \psi ( x_2 , x_3 , t ) , \ {\rm and}
\end{equation}
\begin{equation}
(x_1 , x_2) \in S_t (x_3) \ {\rm if \ and \ only \ if} \
x_1 = \psi ( x_2 , x_3 , t ) \, .
\end{equation}

\noindent
In view of (16) $\cdots$ (19), we have
\begin{equation}
\nu = \left(
1 , \, - \frac{\partial}{\partial x_2} \psi , \,
- \frac{\partial \psi}{\partial x_3}
\right) {\Bigg \slash} \,
\sqrt{1 +
\left( \frac{\partial}{\partial x_2} \, \psi  \right)^2 \, + \,
\left(
\frac{\partial}{\partial x_3} \, \psi \right)^2} \ {\rm and}
\end{equation}

\begin{equation}
\tilde{\nu} \, = \,
\left(
1 , \, - \frac{\partial}{\partial x_2} \, \psi , \,
0 \right) \,
{\Bigg \slash} \,
\sqrt{1 \, + \, \left(
\frac{\partial}{\partial x_2} \, \psi \right)^2} \ {\rm on} \
{\cal{S}} \cap ( U \times J ).
\end{equation}

\noindent
From (10), (17), we have
$\theta ( \psi ( x_2 , x_3 , t ) , x_2 , x_3 , t ) = 0$ on ${\cal{S}} \cap ( U \times J )$.  Differentiating in $t$, we obtain
\[
\left(
\frac{\partial}{\partial t} + \left(
\frac{\partial \psi}{\partial t} \, ( x_2 , x_3 , t ) , 0 , 0 \right) \cdot
\nabla_x \right) \, \theta = 0 \ {\rm on} \
{\cal{S}} \cap ( U \times J ) .
\]

\noindent
Subtracting this from (13), we find that
$\left[ \left(
\frac{\partial \psi}{\partial t} \,
( x_2 , x_3 , t ) , 0 , 0 \right) \, - \sigma \nu \right]$ is orthogonal to
$\nabla_x \theta$, hence also to $\nu$.  This yields the formula

\begin{equation}
\sigma \, = \,
\left( \frac{\partial \psi}{\partial t} \, ( x_2 , x_3 , t ) , 0 , 0 \right)
\, \cdot \, \nu \, = \,
\left( \frac{\partial \psi}{\partial t} \right)
\, {\Bigg \slash}
\,
\sqrt{1 \, + \,
\left( \frac{\partial \psi}{\partial x_2} \right)^2 \, + \,
\left( \frac{\partial \psi}{\partial x_3} \right)^2 } \, .
\end{equation}

\noindent
Similarly, subtracting the two equations (13), we learn that
$( \sigma \nu - \, \tilde{\sigma} \tilde{\nu})$ is orthogonal to $\nabla_x \theta$, hence also to $\nu$.  Therefore,

\begin{equation}
\sigma \, = \,
\tilde{\sigma} \, ( \nu \cdot \tilde{\nu}) \ {\rm on} \
{\cal{S}} \, .
\end{equation}

\noindent
From (20) $\cdots$ (23), we obtain

\begin{equation}
\tilde{\sigma} \, = \,
\left( \frac{\partial \psi}{\partial t} \right) {\Bigg \slash} \,
\sqrt{1 \, + \, \left( \frac{\partial \psi}{\partial x_2} \right)^2} \
{\rm on} \ {\cal{S}} \cap ( U \times J ) \, .
\end{equation}

\noindent
Now we can read off (14) and (15).  In fact, (18) gives

\[
\frac{d}{dt} \left[ \:
\displaystyle{\int\limits_{\Omega_t ( x_3)}} \,
F d \, {\rm (Area)} \right] \, = \, \frac{d}{dt}
\,
\left[
\displaystyle{
\hspace{.35in}
{\int\!\!\!\!\int}_
{_
{_
{\vspace{.80in}
\hspace{-.45in}
{x_1 < \psi ( x_2 , x_3 , t )}
}}}} \,
F dx_1 \, dx_2 \right] \, = \,
\int \, F ( \psi , x_2 , x_3 ) \,
\frac{\partial \psi}{\partial t} \, dx_2 , \]

\noindent
and (19), (24) yield

\[
\displaystyle{\int\limits_{S_t (x_3)} F \tilde{\sigma}} \,
d \ {\rm (length)} \, = \,
\int \, F \tilde{\sigma} \, \cdot \,
\sqrt{1 + \left(
\frac{\partial \psi}{\partial x_2}
\right)^2} \,
dx_2 \, = \,
\int F ( \psi , x_2 , x_3 ) \,
\frac{\partial \psi}{\partial t} \, dx_2 , \]

\noindent
proving (14).  Similarly, (17) gives

\[
\int\limits_{S_t (I)} F \, d \ {\rm ( Area )} \, = \,
\int \, F ( \psi, x_2 , x_3 ) \, \cdot \,
\sqrt{
1 + \left( \frac{\partial}{\partial x_2} \, \psi \right)^2 \, + \,
\left( \frac{\partial}{\partial x_3} \psi \, \right)^2} \, dx_2 \, dx_3 \,
\]

\noindent
while (20) and (21) imply

\[ \nu \cdot \tilde{\nu} \, = \,
\sqrt{1 + \left( \frac{\partial}{\partial x_2} \, \psi \right)^2} \,
{\Bigg\slash} \,
\sqrt{1 + \left( \frac{\partial}{\partial x_2} \, \psi \right)^2 +
\left(
\frac{\partial}{\partial x_3} \, \psi \right)^2}, \]

\noindent
so that (19) yields

\[
\int\limits_{x_3 \in I} \,
\left\{
\int\limits_{S_t ( x_3)} \:
\frac{F}{\nu \cdot \tilde{\nu}} \,
d \ {\rm ( length )} \
\right\}
\, dx_3 \, = \,
\int\limits_{x_3 \in I} \,
\left\{
\int \, \left(
\frac{F}{\nu \cdot \tilde{\nu}} \right) \,
\sqrt{
1 + \left(
\frac{\partial}{\partial x_2} \,
\psi \right)^2} \, dx_2 \right \} \, dx_3 \]

\[
= \, \displaystyle{\int} \,
F ( \psi , x_2 , x_3 ) \, \cdot \,
\sqrt{1 +
\left(
\frac{\partial}{\partial x_2} \, \psi \right)^2  \, + \,
\left(
\frac{\partial}{\partial x_3} \, \psi \right)^2} \, dx_2 \, dx_3 .\]

\noindent
This completes
the proof of (15).  Note that (23) allows us to rewrite (15) in the form

\begin{equation}
\int\limits_{S_t (I)} \, F \, d {\rm (Area )} \, = \,
\int\limits_{x_3 \in I} \,
\left\{ \:
\int\limits_{S_t ( x_3)} \, F \,
\frac{\tilde{\sigma}}{\sigma} \,
d \, {\rm (length)}
\right\} dx_3 \, .
\end{equation}

\section{\sc Proof of the Theorem}

We retain the notation of the previous sections.  We will define a
time-dependent interval

\begin{equation}
J_t \, = \, [ A ( t ) , B ( t ) ] \subset I_3
\end{equation}

\noindent
and establish an obvious formula for the time derivative of
Vol $\Omega_t ( J_ t)$.  We assume that the endpoints
$A ( t )$, $B(t)$ are $C^1$ functions of $t$.  We have

\[ {\rm Vol} \,\Omega_t (J_t) \, = \,
\int\limits_{x_3 \in J_t} \, {\rm Area} \,
\Omega_t ( x_3 ) dx_3 , {\rm so \ that} \]

\[
\frac{d}{dt} \, {\rm Vol} \, \Omega_t ( J_t ) \, = \,
B^\prime ( t ) \, {\rm Area} \, \Omega_t ( B(t)) - A^\prime ( t ) \,
{\rm Area} \, \Omega_t (A ( t )) \, + \,
\int\limits_{x_3 \in J_t} \,
\frac{\partial}{\partial t} \,
{\rm Area} \,
\Omega_t ( x_3) \, dx_3 \, .  \]

\noindent
Applying (14) with $F \equiv 1$, we find that

\[
\frac{d}{dt}  {\rm Vol} \,
\Omega_t ( J_t )  =
B^\prime (t) \, {\rm Area} \,
\Omega_t ( B ( t )) \,
 - A^\prime(t) \,
{\rm Area} \, \Omega_t ( A ( t ))  +
\int\limits_{x_3 \in J_t}
\left\{
\int\limits_{S_t ( x_3)}
\tilde{\sigma} \, d {\rm ( length )} \right\} \, dx_3 \, .
\]

\noindent
In view of (25) (with $F \equiv \sigma$ on ${\cal{S}}$), this is
equivalent to
\begin{equation}
\frac{d}{dt} \,
{\rm Vol} \, \Omega_t ( J_t ) =
B^\prime(t) \,
{\rm Area} \,
\Omega_t (B(t)) \,
- A^\prime ( t ) \,
{\rm Area} \,
\Omega_t (A ( t )) +
\int\limits_{S_t (J_t)} \,
\sigma \, d {\rm ( Area )}  \, .
\end{equation}

\noindent
Now we bring in the hypothesis that $\Omega_t$ moves with a divergence-free
$C^1$ velocity field $u$.  From (6) and (13), we see that
$( \sigma \nu - u ) \cdot \nabla_x \theta = 0$ on $S_t (J_t)$.  Thus
$( \sigma \nu - u )$ is orthogonal to $\nu$, so that
$\sigma = u \cdot \nu$ on $S_t ( J_t)$, and (27) may be rewritten as
\begin{equation}
\frac{d}{dt} \, {\rm Vol} \,
\Omega_t (J_t) =
B^\prime ( t ) \, {\rm Area} \,
\Omega_t (B(t)) \,
- A^\prime (t) \, {\rm Area} \, \Omega_t (A(t)) +
\int\limits_{S_t(J_t)} \, u \cdot \nu d \, {\rm (Area)} \, .
\end{equation}

\noindent
On the other hand, since $u = (u_1 , u_2 , u_3)$ is divergence-free, the
divergence theorem yields

\[
0 = \int\limits_{\Omega_t ( J_t)} \,
( \nabla_x \cdot u) \,
d ({\rm Vol}) =
\int\limits_{S_t(J_t)} \,
u \cdot \nu \, d {\rm (Area)} \,
+
\int\limits_{\Omega_t (B(t))} \, u_3 \, d {(\rm Area)} \,
- \int\limits_{\Omega_t (A ( t))} \, u_3 d \, {\rm(Area)} \, .
\]

\noindent
Hence, (28) may be rewritten in the form
\begin{equation}
\frac{d}{dt} \, {\rm Vol} \, \Omega_t (J_t) =
\int\limits_{\Omega_t (B(t))} \,
[ B^\prime ( t ) - u_3 ( x , t ) ] \, d \,
{\rm (Area)} \,
- \int\limits_{\Omega_t ( A ( t ) )} \,
[ A^\prime ( t ) - u_3 ( x , t )] \, d {\rm (Area)} \, .
\end{equation}

\noindent
This is our final formula for the time derivative of Vol $\Omega_t (J_t)$.
It is intuitively clear.

\noindent
We now pick the time-dependent interval
$J_t = [ A ( t ) , B ( t ) ] \subset I_3$. Let $I_3 = [ a , b ]$, and let
$t_0 \in ( 0 , T)$ be a time to be picked below.  We define

\begin{equation}
B ( t ) = b - \int\limits^{T}_{t} \, \max_{x \in Q} \, | u ( x , \tau ) | \, d \tau \, ,
\end{equation}
\noindent
and
\begin{equation}
A ( t ) = a + \int\limits^T_t \, \max_{x \in Q} \, | u ( x , \tau ) | \, d \tau \, .
\end{equation}

\noindent
We are assuming that $u ( x , \tau )$ is continuous on $Q \times [ 0 , T )$,
and that
$\int\limits_{0}^{T} \max\limits_{x \in Q}
| u ( x , \tau ) | d \tau < \infty$.
It follows that $A ( t ), B ( t )$ are $C^1$ functions on $[0 , T)$, and that

\begin{equation}
a \leq A(t) < B ( t ) \leq  b \ {\rm for} \ t \in [ t_0 , T ) \, ,
\end{equation}

\noindent
provided we pick $t_0$ close enough to $T$.  We pick $t_0$ so that (32) holds.
Thus, $\Omega_t (J_t) \subset Q$ for $t \in [t_0 , T)$.  Immediately from (30), (31), we obtain

\begin{equation}
B^\prime (t) = \,
-A^\prime (t) = \max\limits_{x \in Q} \,
|u ( x , t ) | \, \geq \,
\max\limits_{x \in \Omega_t (A ( t )) \, \cup \, \Omega_t (B ( t ))} \: | u_3 ( x , t ) |
\end{equation}

\noindent
(recall $u = (u_1 , u_2 , u_3 ))$.

\noindent
From (29) and (33) we see at once that

\begin{equation}
\frac{d}{dt} \, {\rm Vol} \, \Omega_t ( J_t ) \, \geq \,
0 \ {\rm for } \ t \in [ t_0 , T ) \, .
\end{equation}

\noindent
On the other hand, (4) and (32) show that
${\rm Vol} \, \Omega_{t_0} ( J_{t_0} ) \, > \, 0$.

\noindent
Consequently,

\[
\liminf\limits_{t \to T-} \, {\rm Vol} \, \Omega_t \, \geq \,
\liminf\limits_{t \to T-} \,
{\rm Vol} \, \Omega_t ( J_t ) \, \geq \, {\rm Vol} \, \Omega_{t_0} ( J_{t_0} ) \, > \,0 \, .\]

\noindent
The proof of our theorem is complete. \hfill $\blacksquare$

\vfill
\newpage
\section{ \sc References}
\begin{enumerate}
\item J. Beale, T. Kato, and A. Majda, ``Remarks on the breakdown of smooth solutions for the 3D Euler equations,'' {\it Comm. Math. Phys.}, {\bf 94}, pp.61-64, (1984).
\item P. Constantin, C. Fefferman, and A. Majda, ``Geometric constraints
on potentially singular solutions for the 3D Euler equations,''
{\it Commun. Part. Diff. Eq.}, {\bf 21}, pp.559-571, (1996).
\item D. Cordoba and C. Fefferman, ``Scalars convected by a 2D
incompressible flow,'' (to appear).
\item D. Cordoba and C. Fefferman, ``Behavior of several 2D fluid equations in
singular scenarios,'' submitted to {\it Proc. Nat. Acad. Sci.}, U.S.A.
\end{enumerate}
\vfill
\today
\end{document}